\documentclass[10pt]{article}

\usepackage{amsmath,amssymb}
\newtheorem{thm}{Theorem}
\newtheorem{prop}{Proposition}
\newtheorem{lem}{Lemma}
\newtheorem{coro}{Corollary}
 
\newtheorem{conj}{Conjecture}
\usepackage{mathtools}

\DeclarePairedDelimiter\floor{\lfloor}{\rfloor}

\def\bbr{{\mathbb R}}

\def\bbe{{\mathbb E}}
\def\bbp{{\mathbb P}}

\def\itm#1{\par\parindent=20pt\noindent
   \hangindent\parindent\hbox to\parindent{#1\hss}\ignorespaces}

\def\dperp{\mathrel{\perp\kern-1em\perp}}

\author{Philippe Marchal}
\title
{Rectangular Young tableaux and the Jacobi ensemble}

\begin{document}
\maketitle
\begin{abstract}

It has been shown by Pittel and Romik that the random surface associated with a 
large rectangular Young tableau converges to a deterministic limit. 
We study the fluctuations from this limit along the edges of the rectangle.
We show that in the corner,  these fluctuations are gaussian wheras, away 
from the corner and when the rectangle is a square, the fluctuations
are given by the Tracy-Widom distribution. Our method 
is based on a connection with the Jacobi ensemble.\\

\end{abstract}


\section{Statement of the results}

From a formal point of view, a rectangular Young tableau of size $(m, n)$
can be defined as an $mn$-tuple
of integers $(X_{1,1},\ldots X_{m,n})$ 
satisfying, for all $i,j$,
$$
X_{i,j}< \min(X_{i,j+1}, X_{i+1,j})
$$
$$\{X_{1,1},\ldots X_{m,n}\}=\{1,2,\ldots mn\}$$ 
We denote by $\mathcal{X}_{m,n}$ the set of $mn$-tuples of this form.

Pittel and Romik \cite{pitt} 
studied the limit shape of the surface associated with 
a random large rectangular Young tableau. Fix a real $t>0$,
and consider a Young tableau of size $(n, \floor{tn})$, chosen 
uniformly at random. Then they proved that for all reals 
$(r, s)\in [0,1]^2$, $X_{\floor{rn}, \floor{stn}}/tn^2$ 
converges in law as $n\to \infty$ to a 
deterministic quantity $g(r,s,t)$, which is 
the solution of a minimization problem.
We want to study the fluctuations of this shape along the edges of the 
rectangle. First, we state a completely explicit result on the corner:

\begin{thm}
Let $(X_{1,1},\ldots,X_{m,n})$
be a uniform random variable on $\mathcal{X}_{m,n}$.
Then for $n\leq k\leq mn-m+1$,
\begin{equation}\label{t1}
\bbp(X_{1,n}=k)=\frac{\binom{k-1}{n-1}\binom{mn-k}{m-1}}
{\binom{mn}{m+n-1}}
\end{equation}
\end{thm}

Let us compare this result with the following model. Let $R$ be a 
rectangle of size $((m-1)(n-1), m+n-1)$ on an integer lattice. 
Consider the set of paths going from the south-west corner
to the north-east corner
with only north and east steps. Choose such a path uniformly at random.
Then  it is an elementary exercise to check that the probability for the 
$n$-th north step to be the $k$-th step is equal to \eqref{t1}.
However, we have not been able to find a combinatorial link between this simple 
model and the corner of a rectangular Young tableau.

A consequence of Theorem 1 is the following:

\begin{coro}
Fix some real $t\in \bbr_+^*$. For every $n$, put $m_n=\floor{tn}$ and 
let $(X^{(n)}_{1,1},\ldots,X^{(n)}_{m_n,n})$
be a uniform random variable on $\mathcal{X}_{m_n,n}$.  Then as $n\to\infty$,
$$\frac{\sqrt{2}(1+t)\left(X^{(n)}_{1,n}-\bbe X^{(n)}_{1,n}\right)}{n^{3/2}}
\stackrel{law}{\to}G
$$ 
where $G$ is a standard gaussian random variable.
\end{coro}

Our second result only holds 
when the rectangle is a square and deals with the fluctuations on the edges.

\begin{thm}
For each $n$, let $(X^{(n)}_{1,1},\ldots,X^{(n)}_{n,n})$
be a uniform random variable on $\mathcal{X}_{n,n}$. 
Then there exists a function $r:(0,1)\to \bbr_+^*$
such that for $t\in(0,1)$, as $n\to\infty$,
$$\frac{r(t)\left(X^{(n)}_{\floor{tn},n} -\bbe X^{(n)}_{\floor{tn},n}\right)}{n^{4/3}}
\stackrel{law}{\to} TW
$$ 
where $TW$ has the Tracy-Widom distribution.
\end{thm}

Pittel-Romik's result tells us that
$n^{-2}\bbe X^{(n)}_{\floor{tn},n}$
converges to the value of the limit shape at $(t,1)$, 
namely $(1+\sqrt{2t-t^2})/2$.

Our method consists in studying Young tableaux in a slightly modified framework.
As the asymptotic shape only makes sense when the ``height'' is renormalized,
that is, the $X_{i,j}$ are divided by $mn$, 
it seems natural to work directly in
a continuous framework.
The formal setup is the following. \\

\noindent
{\bf Definition}. For a pair of integers 
$(m, n)$, let $\mathcal{Y}_{m,n}$ be the set of $mn$-tuples  
$(Y_{1,1},\ldots Y_{m,n})$ of reals in
$[0,1]$ satisfying, for all $i,j$,
\begin{equation}\label{e1}
Y_{i,j}< \min(Y_{i,j+1}, Y_{i+1,j})
\end{equation}
We want to study uniform random variables on $\mathcal{Y}_{m,n}$. 
Denote by  $\Delta(x_1,\ldots x_k)$
the Vandermonde of the $k$-tuple $ (x_1,\ldots x_k)$:
$$
\Delta(x_1,\ldots x_k)=\prod_{1\leq i<j\leq k}(x_i-x_j)
$$
When the rectangle is a square, we get:

\begin{thm}
Let $n$ be a positive integer and $(Y_{1,1},\ldots,Y_{n,n})$ be a 
uniform random variable on the set $\mathcal{Y}_{n,n}$.
Then for every $k\in[1, n]$, the $k-tuple$
$$
(Y_{1,n-k+1},Y_{2,n-k+2},\ldots,Y_{k-1,n-1},Y_{k,n})
$$
has a marginal density proportional to
$$
{\bf 1}_{\{x_1\leq x_2\ldots \leq x_k \}}
\Delta(x_1,\ldots x_k)^2\prod_{i=1}^k x_i^{n-k}(1-x_i)^{n-k}
$$
\end{thm}

For a general rectangle, the result reads:

\begin{thm}
Let $m < n$ be two positive integers. Let $(Y_{1,1},\ldots,Y_{m,n})$ be a 
uniform random variable on the set $\mathcal{Y}_{m,n}$.

(i)  Let $k\in[1, m]$ be
an integer. Then the $k$-tuple
$$
(Y_{1,n-k+1},Y_{2,n-k+2},\ldots,Y_{k-1,n-1},Y_{k,n})
$$
has a marginal density proportional to
$$
{\bf 1}_{\{x_1\leq x_2\ldots \leq x_k \}}
\Delta(x_1,\ldots x_k)^2\prod_{i=1}^k x_i^{n-k}(1-x_i)^{m-k}
$$

(ii)
If  $m<k\leq n$, then the $m-tuple$
$$
(Y_{1,n-k+1},Y_{2,n-k+2},\ldots,Y_{m-1,n-k+m-1},Y_{m,n-k+m})
$$
has a marginal density proportional to
$$
{\bf 1}_{\{x_1\leq x_2\ldots \leq x_k \}}
\Delta(x_1,\ldots x_k)^2\prod_{i=1}^k x_i^{n-k}(1-x_i)^{k-m}
$$
\end{thm}

Of course, for a diagonal of the form $(Y_{m-k+1,1},Y_{m-k+2,2},\ldots,Y_{m,k})$, 
we get a similar expression as in (i).
The densities appearing in these two theorems belong to the general class 
called the 
Jacobi ensemble. This ensemble also appears in various models,
among others the MANOVA procedure in statistics \cite{muir}, 
log-gas theory \cite{forr}, Wishart matrices and random 
projections \cite{coll}. For a detailed account on random matrices, we refer 
to \cite{and}. 

Using Theorems 3 and 4 together with known results on the Jacobi ensemble
enables us to derive the results stated above. Moreover, the deterministic 
limit shape can also be recovered this way, 
see Section 4. An alternative form of Theorem 2 in the continuous setting
is as follows:

\begin{coro}
For each integer $n$, let $(Y^{(n)}_{1,1},\ldots,Y^{(n)}_{n,n})$
be a uniform random variable on $\mathcal{Y}_{n,n}$. 
Then for $t\in(0,1)$, with the same function $r$ as in Theorem 2, as 
$n\to\infty$,
$$\frac{r(t)(Y^{(n)}_{\floor{tn},n} -\bbe Y^{(n)}_{\floor{tn},n})}{n^{4/3}}
\stackrel{law}{\to} TW$$ 
where TW has the
Tracy-Widom distribution. 
\end{coro}

Remark that $\mathcal{Y}_{n,n}$ is a convex subset (indeed, a compact polytope)
of $\bbr^{n^2}$. From this point of view, Corollary 2 can be seen as a result on
the projection of the uniform measure on a convex set in high dimension. This
is reminiscent of the classical result saying that if $(X_1,\ldots X_n)$ is 
a random
vector distributed according to the uniform measure on the euclidean 
$n$-dimensional sphere with radius $\sqrt{n}$, then $X_1$ is asymptotically 
gaussian.

 The remainder of this note is organized as follows.
We prove Theorems 3 and 4 in Section 2. Theorem 1 is derived in Section 3. 
Section 4 is devoted to the proofs of the asymptotic results. Some concluding 
remarks are made in Section 5

\section{Diagonals of continuous tableaux}

We prove here Theorem 3. The proof of Theorem 4 is the same and is omitted.
The basic idea is to use a random generation algorithm filling the tableau 
diagonal by diagonal, using conditional densities.

We begin with a preliminary lemma. 
Let $n\geq 2$  and define by induction the following polynomials: 
$g_1(x)=1
$
and for $i\leq n-1$,
$$
g_{i+1}(x_1,\ldots x_{i+1})=\int_{x_1}^{x_2}dy_1\ldots
\int_{x_i}^{x_{i+1}}dy_n\ g_i(y_1,\ldots, y_i)
$$
while for $i\in[n, 2n-2]$,
$$
g_{i+1}(x_1,\ldots x_{2n-i-1})=\int_0^{x_1}dy_1\int_{x_1}^{x_2}dy_2\ldots
\int_{x_{2n-i-1}}^1dy_{2n-i}\ 
g_i(y_1,\ldots, y_{2n-i})
$$
 
\begin{lem}
(i) For every $i\in [1,n]$, there exists a constant $c_i$ such that
$$
g_i(x_1,\ldots x_{i})= c_i\Delta(x_1,\ldots x_i)
$$

(ii) For $i\in[n, 2n-1]$, there exists a constant $c_i$ such that
$$
g_{i}(x_1,\ldots x_{2n-i})=c_i\Delta(x_1,\ldots x_{2n-i})\prod_{j=1}^{2n-i}
x_j^{i-n}(1-x_j)^ {i-n}
$$
\end{lem}

{\bf Proof}

An elementary proof of (i) can be found in Baryshnikov \cite{bary}.
To deduce (ii), we proceed by induction. For $m\leq n$, define 
$K_m(\varepsilon,y_1,\ldots y_{n-m}) $ as the integral
\begin{eqnarray*}
&&\int_{-\varepsilon}^0 dr_1 \ldots\int_{-m\varepsilon}
^{-(m-1)\varepsilon}dr_m
\int_1^{1+\varepsilon}ds_1\ldots  \int_{1+(m-1)\varepsilon}^{1+m\varepsilon}ds_m\\
&&\ \ \ 
\Delta(r_1,\ldots r_m, y_1,\ldots y_{n-m},s_1,\ldots  s_m)
\end{eqnarray*}
We want to evaluate
$$
I:=\int_0^{x_1}dy_1\ldots \int_{x_{n-m-1}}^1dy_{n-m} 
K_m(\varepsilon,y_1,\ldots y_{n-m})
$$
On the one hand, using (i) easily gives
\begin{equation}\label{I}
I\sim \varepsilon^{-m(m+1)}C(m,n) \Delta( x_1,\ldots,x_{n-m-1})
\prod_{i=1}^{n-m-1} x_i^{m+1}(1-x_i)^{m+1}
\end{equation}
for some positive constant $C(m,n)$, where the equivalent is 
uniform over all $(n-m)$-tuples ($x_1, \ldots x_{n-m})$
satisfying
\begin{equation}\label{unif}
\sqrt{\varepsilon}\leq x_1 \leq \ldots\leq x_{n-m}\leq 1-\sqrt{\varepsilon}
\end{equation}
On the other hand,
$$
\Delta(r_1,\ldots y_1,\ldots  s_m)
\sim
\Delta(y_1, \ldots y_{n-m}) \Delta(r_1,\ldots  s_1,\ldots  s_m)
\prod_{i=1}^{n-m} y_i^m(1-y_i)^m
$$
where the equivalent is uniform over all $(n+m)$-tuples 
$(r_1,\ldots  y_1,\ldots s_m)$
satisfying
$$
\sqrt{\varepsilon}\leq y_1 \leq \ldots\leq y_{n-m}\leq 1-\sqrt{\varepsilon}
$$
$$-\varepsilon<r_1<0,\ldots,-m\varepsilon< r_m< -(m-1)\varepsilon,  
1<s_1<1+\varepsilon,\ldots $$

As a consequence,
$$
K_m(\varepsilon,y_1,\ldots y_{n-m})\sim C'(m,n)
\varepsilon^{m(m+1)}\Delta(y_1, \ldots, y_{n-m}) 
\prod_{i=1}^{n-m} y_i^m(1-y_i)^m
$$
where $C'(m,n) $ is the constant such that
$$
 \int_{-\varepsilon}^0 dr_1 \ldots
\int_{1+(m-1)\varepsilon}^{1+m\varepsilon}ds_m
\Delta(r_1,\ldots  s_1,\ldots  s_m)
\sim  C'(m,n)\varepsilon^{m(m+1)}
$$
It follows that uniformly over all $(n-m)$-tuples ($x_1, \ldots x_{n-m})$
satisfying \eqref{unif},
$$I \sim C'(m,n)
\varepsilon^{m(m+1)}\int_0^{x_1}dy_1\ldots \int_{x_{n-m-1}}^1dy_{n-m}
 \Delta( y_1,\ldots,y_{n-m})\prod_{i=1}^{n-m} y_i^m(1-y_i)^{m}
$$
Comparing this estimate with \eqref{I} yields (ii).\hfill $\Box$

We now describe an algorithm generating random elements in 
$\mathcal{Y}_{n,n}$. A similar algorithm can be used to generate
random permutations with a prescribed profile of ascents and descents
\cite{mar}. In the remainder of this section, $k$ and $n$ are the 
integers in the statement of  Theorem 3.
For $i\in[1,n]$, we denote 
$$
D_i= (Y_{1,n-i+1},Y_{2,n-i+2},\ldots,Y_{i,n})
$$
while for $i\in[n+1,2n-1]$,
$$
D_{i}= (Y_{i-n+1,1},Y_{i-n+2,2},\ldots,Y_{n,2n-i})
$$
If $(x_1, \ldots x_j)$ and 
$(y_1, \ldots y_{j+1})$
are two sequences of reals, we say that they are interlacing in $[0,1]$ if
$
0\leq y_1\leq x_1\leq y_2\ldots \leq x_j\leq  y_{j+1}\leq 1
$.
We denote the event that this interlacing relation is satisfied by
$$Inter((x_1, \ldots x_j),(y_1, \ldots y_{j+1}))$$

{\bf Algorithm}
\begin{itemize}
\item
Choose the diagonal $D_k$
at random according to the density
$$
\frac{{\bf 1}_{\{0\leq x_1\leq x_2\ldots \leq x_k \leq 1\}}
g_k(x_1,\ldots x_k)g_{2n-k}(x_1,\ldots x_k)}{Z_k}
$$ where 
$$
Z_k=
{\int_0^1 dx_k \int_0^{x_k} dx_{k-1}\ldots \int_0^{x_2} dx_1
g_k(x_1,\ldots x_k)g_{2n-k}(x_1,\ldots x_k)}
$$
\item
By induction, for $i$ from $k$ down to $2$, conditional on 
$D_i$, choose $D_{i-1}$
according to the conditional density
\begin{equation}\label{g}
\frac{g_{i-1}(x_1,\ldots, x_{i-1}){\bf 1}_{Inter(D_i, (x_1,\ldots, x_{i-1}))}}
{g_{i}(D_i)}
\end{equation}
\item
By induction, for $i$ from $k$ to $n-1$, conditional on 
$D_i$, choose $D_{i+1}$
according to the conditional density
$$
\frac{g_{i+1}(x_1,\ldots, x_{i+1})
{\bf 1}_{Inter(D_i, (x_1,\ldots, x_{i+1}))}}
{g_{i}(D_i) }
$$
\item
By induction, for $i$ from $n$ to $2n-2$, conditional on 
$D_i$, choose $D_{i+1}$
according to the conditional density
$$
\frac{g_{i+1}(x_1,\ldots, x_{2n-i-1})
{\bf 1}_{Inter(D_i, (x_1,\ldots, x_{2n-i-1}))}}
{g_{i}(D_i) }
$$
\end{itemize}

First, remark that the conditional densities used by the algorithm are indeed
probability densities. That is, they are measurable, positive functions
and their integral is 1. The latter fact is easy to verify: for instance, by 
definition of $g_i$,
$$
g_i(D_i)=g_{i}(Y_1,\ldots Y_{i})=\int_{Y_1}^{Y_2}dx_1\ldots
\int_{Y_{i-1}}^{Y_{i}}dx_{i-1}\ g_{i-1}(x_1,\ldots, x_{i-1})
$$
and thus, using \eqref{g}, we get
$$
\int_{[0,1]^{i-1}}dx_1\ldots dx_{i-1}\frac{g_{i-1}(x_1,\ldots, x_{i-1})
{\bf 1}_{Inter(D_i, (x_1,\ldots, x_{i-1}))}}
{g_i(D_i)}=1
$$
We claim that the algorithm yields a random element of  
$\mathcal{Y}_{n,n}$ with the uniform measure. Indeed, by construction,
the $n^2$-tuple generated by the algorithm has a density which is the 
product of the conditional densities of the diagonals 
$D_1, D_2\ldots D_{2n-1}$. Hence this 
density is given by

\begin{eqnarray*}
&&\frac{{\bf 1}_{\{0\leq x_1\leq \ldots \leq x_k \leq 1\}}g_k(D_k)g_{2n-k}(D_k)}{Z_k}\\
&&\times \prod_{i=2}^{k}\frac{g_{i-1}(D_{i-1})
{\bf 1}_{Inter(D_i, D_{i-1})}}{g_i(D_i)}\prod_{i=k+1}^{2n-1}\frac{g_{i}(D_{i})
{\bf 1}_{Inter(D_i, D_{i-1})}}{g_{i-1}(D_{i-1})}
\end{eqnarray*}
The expression above is a telescopic product and after simplification, 
we find that the density  is constant on the set 
$\mathcal{Y}_{n,n}$. This proves our claim. 

Finally, the density of  $D_k$ is proportional to
\begin{eqnarray*}
&&{\bf 1}_{\{x_1\leq \ldots \leq x_k \}}
g_k(x_1,\ldots x_k)g_{2n-k}(x_1,\ldots x_k) \\
&& ={\bf 1}_{\{x_1\leq \ldots \leq x_k \}}
\Delta(x_1,\ldots x_k)^2\prod_{i=1}^k x_i^{n-k}(1-x_i)^{n-k}
\end{eqnarray*}
according to Lemma 1. This proves Theorem 3.\\

\section{The law of the corner}

One can relate the discrete and the continuous model of Young tableaux. 
To construct a continuous Young tableau $(Y_{i,j})$ of size $(m,n)$ 
from a discrete Young tableau $(X_{i,j})$ of the same size, proceed as follows:

\begin{itemize}
\item
Let $(X_{i,j})$ be a uniform random variable on $\mathcal{X}_{m,n}$.
\item
Let 
$(Z_1\leq \ldots\leq Z_{mn})
$
be the increasing reordering of $mn$ independent, uniform 
random variables on $[0,1]$, independent of $(X_{i,j})$.
\item
For every pair $(i,j)$, let $k(i,j)$ be the integer satisfying $X_{i,j}=k(i,j)$.
Then put $Y_{i,j}=Z_{k(i,j)}$.
\end{itemize}

\begin{prop}
Consider the $(mn)$-tuple $(Y_{i,j})$ constructed as above. Then

(i)
$(Y_{i,j})$ is 
distributed according to the uniform measure on $ \mathcal{Y}_{m,n}$,

(ii) For every 
pair $(i,j)$,

$$
\bbe\left(Y_{i,j}-\frac{X_{i,j}}{mn+1}\right)^2\leq  \frac{1}{mn+1}
$$

\end{prop}

The proof of (i) is elementary and (ii) follows from a simple variance 
computation, using the fact that the density of $Z_k$ is
$$
h_k(x)=x^{k-1}(1-x)^{mn-k}\frac{(mn)!}{(k-1)!(mn-k)!}
$$

As a consequence of Proposition 1, let $f_{i,j}$ be the marginal 
density of $Y_{i,j}$.
For every $1\leq k\leq mn$  put
$p_{i,j}(k)=\bbp(X_{i,j}=k)
$. 
Using Proposition 1, we get that the 
density  $f_{i,j}$ is equal to
$$
f_{i,j}(x)=
\sum_{k=1}^{mn}p_{i,j}(k)h_{k}(x)
$$
where $h_k(x)$ is the density of $Z_k$.
Thus
\begin{equation}\label{decomp}
f_{i,j}(x)=
\sum_{k=1}^{mn}p_{i,j}(k)x^{k-1}(1-x)^{mn-k}\frac{(mn)!}{(k-1)!(mn-k)!}
\end{equation}

This way one can deduce the probabilities in the discrete
model from the densities in the continuous model.  For the case
$i=1$, $j=n$, according to Theorem 4 (i), 
$$
f_{1,n}(x)=x^{n-1}(1-x)^{m-1}\frac{(m+n-1)!}{(m-1)!(n-1)!}
$$
To obtain the desired decomposition,  divide both sides of
\eqref{decomp}  by $(1-x)^{mn-1}$ and use the change of variables 
$y=x/(1-x)$ to get
$$
\frac{(m+n-1)!}{(m-1)!(n-1)!}y^{n-1}(1+y)^{mn-n-m+1} =
\sum_{k=1}^{mn}p_{1,n}(k)y^{k-1}\frac{(mn)!}{(k-1)!(mn-k)!}
$$
Identifying the coefficient of $y^{k-1}$, 
we find that if $n\leq k\leq mn-m+1$, then
$$
p_{1,n}(k)=c_{m,n}
\frac{(k-1)!(mn-k)!}{(k-n)!(mn-m-k+1)!}
$$
with
$$
c_{m,n}=
\frac{(m+n-1)!(mn-m-n+1)!}
{(mn)!(m-1)!(n-1)!}
$$
and $p_{i,j}(k)=0$ otherwise. This proves Theorem 1.

\section{Asymptotic results}

\subsection{Proof of Theorem 2} 

For every $n$, let $(Y^{(n)}_{1,1},\ldots, Y^{(n)}_{n,n})$ be distributed
according to the uniform measure on  
$\mathcal{Y}_{n,n}$. Collins' results on the ``soft edge'' of the Jacobi 
ensemble \cite{coll}
can be translated in our context as follows. For every $t\in(0,1)$, 
there exists a sequence $(s_n(t))$ and a constant $r(t)$ such that, for every 
$k$, 
as $n$ goes to infinity,
\begin{equation}\label{airy}
r(t)n^{2/3}\left(Y^{(n)}_{\floor{tn},n}-s_n(t), \ldots, 
Y^{(n)}_{\floor{tn}-k+1,n-k+1}-s_n(t)\right)
\end{equation}
converges in distribution to  a $k$-tuple 
$
(A_1,\ldots,A_k)
$
which is the truncation of a random point process 
$$(A_1> A_2>\ldots >A_k>\ldots)$$ 
distributed 
according to the Airy ensemble.
Moreover, the estimate (ii) of Proposition 1 entails that \eqref{airy} and
$$
r(t)n^{2/3}\left(\frac{X^{(n)}_{\floor{tn},n}}{n^2}-s_n(t), \ldots, 
\frac{X^{(n)}_{\floor{tn}-k+1,n-k+1}}{n^2}-s_n(t)\right)
$$
have the same limit. This yields Theorem 2. Remark that Theorem 2 is a 
specialization to the case $k=1$ but that we have in fact a 
multidimensional convergence:

\begin{thm}
With the same assumptions and notation as in Theorem 2, for every $k\geq 1$,
the $k$-tuple
$$
r(t)n^{2/3}\left(\frac{X^{(n)}_{\floor{tn},n}-\bbe X^{(n)}_{\floor{tn},n}}{n^2}, \ldots, 
\frac{X^{(n)}_{\floor{tn}-k+1,n-k+1}-\bbe X^{(n)}_{\floor{tn}-k+1,n-k+1}}{n^2}\right)
$$
converges in distribution towards the largest $k$ values of the Airy process.
\end{thm}

In fact, \cite{coll} gives a general result which applies for every rectangle.
More precisely, for all $m,n,k$, the diagonal
$(Y_{k,n},\ldots, Y_{1,n-k-1})
$
of $\mathcal{Y}_{m,n}$
is a determinantal point process for which \cite{coll} provides estimates both 
on the edge
and in the bulk of the spectrum. However, the convergence towards the Airy 
ensemble 
is only proven in the square case. The extension to the rectangular case would
require some concentration inequalities that do not seem to be available yet.

\subsection{Proof of Corollary 1}
To prove Corollary 1, remark again that  Proposition 1 (ii) enables us to 
deduce 
the limit law of $X^{(n)}_{1,n}$ from the limit law of $Y^{(n)}_{1,n}$. 
The density of $Y^{(n)}_{1,n}$ is
$x^{n-1}(1-x)^{m_n-1}$. Put $r=1/(1+t)$ and
$T_n=\sqrt{n-1}(Y^{(n)}_{1,n}-r)
$.
Then the density of $T_n$ has the form
$$
c_n\left(1+\frac{x}{r\sqrt{n-1}}\right)^{n-1}
\left(1-\frac{x}{(1-r)\sqrt{n-1}}\right)^{\floor{tn}-1}
$$
for some constant $c_n$. Asymptotic estimations easily yield Corollary 1.

\subsection{The deterministic limit shape}

We focus here on the square case but the case of a rectangle can be dealt with
similarly.
Fix a real $t\in (0,1)$, 
let $m(t)=\floor{tn}$ and consider the diagonal 
$$
(D^{(t)}_1, D^{(t)}_2,\ldots D^{(t)}_{m(t)})
:=(Y_{1,n-m(t)+1},Y_{2,n-m(t)+2},\ldots,Y_{m(t),n})
$$
Consider the empirical measure
$$
\mu_n(t)=\frac{1}{m(t)}\sum_{i=1}^{m(t)} \delta_{D^{(t)}_i}
$$
General results on Jacobi ensembles apply in this case and 
we get that $\mu_n(t)$  converges in distribution, as $n$ goes to infinity, to
the deterministic probability measure with density
\begin{equation}\label{a}
f_t(x):=\frac{\sqrt{(\lambda_+(t)-x)(x-\lambda_-(t))}}{\pi x(1-x)}
{\bf 1}_{\{x\in[\lambda_-(t),\lambda_+(t)]\}}
\end{equation}
where
\begin{equation}\label{b}
\lambda_\pm(t)=\frac{1\pm\sqrt{t(2-t)}}{2}
\end{equation}
See for instance the first proposition in \cite{dum}. 
Reformulating this result, we get the following.
Let $r, s\in [0,1]^2$, $r\geq s$ and put $t=1-s+r$. Then as $n$ tends to 
infinity, $X_{\floor{rn}, \floor{sn}}$ converges in law to the Dirac point mass 
$\delta_{g(r,s)}$
where 
$$
g(r,s)= F_t^{-1}(s/t)
$$
the function $F_t^{-1} $ being the inverse of the function
$$
F_t(x)=\int_{\lambda_-(t)}^xf_t(y) dy
$$
and $f_t$ being given by \eqref{a} and  \eqref{b}.
The function $g$ is an alternative fomulation of the limit shape found by 
Pittel and Romik.

Remark however that this result is weaker than Pittel-Romik's, since it only 
gives 
the convergence along a diagonal, whereas Pittel and Romik show a uniform 
convergence
on the whole rectangle. In fact, a major weakness of our method is that it only
allows us to work on a single diagonal, and not on several diagonals 
simultaneously.

\section{Concluding remarks}

Consider the case when the rectangle is a square.
It would be interesting to study the transition between the deterministic 
regime of $X_{n,n}$ and the fluctuations of order $n^{4/3}$ for 
$X_{\floor{tn},n}$, as well as the
transition between the fluctuations of order $n^{4/3}$ for $X_{\floor{tn},n}$ and
the fluctuations of order $n^{3/2}$ for $X_{1,n}$. 
A natural conjecture is the following:

\begin{conj}
Let $(a_n)$ be a nondecreasing sequence with $a_n\to\infty$ as $n\to\infty$.
Then, up to a multiplicative constant,
$$\frac{a_n^{1/6}(X^{(n)}_{\floor{a_n},n} -\bbe X^{(n)}_{\floor{a_n},n})}{n^{3/2}}
\stackrel{law}{\to} TW
$$ 
$$
\frac{X^{(n)}_{\floor{a_n},1} -\bbe X^{(n)}_{\floor{a_n},1}}{a_n^{4/3}}
\stackrel{law}{\to} TW
$$
where $TW$ has the Tracy-Widom distribution.
\end{conj}

Let us explain where this conjecture comes from. The function
$r$ from Theorem 2 can be computed using \cite{coll}:
$$
r(t)=\frac{\sqrt{2}(t(2-t))^{1/6}}{(3-2t+t^2-\sqrt{2t-t^2})^{2/3}}
$$
Computing the asymptotics when $t\to 1$ leads to the second part of the  conjecture.
The first part comes from a link with the random matrix model known as the GUE.
For a fixed $k$,  define a family of variables $(T^{(n)}_{i,k}, 1\leq i\leq k)$ 
by the formula
$$Y_{i, n-k+i}=\frac{1}{2}\left(1-\frac{T^{(n)}_{i,k}}{\sqrt{n}}\right) $$
Then it follows immediately from Theorem 3 that, as $n$ goes to infinity, 
the renormalized diagonal
$(T^{(n)}_{i,k},\ldots T^{(n)}_{k,k})$ has a limit density proportional to
$$
{\bf 1}_{\{x_1\leq x_2\ldots \leq x_k \}}
\Delta(x_1,\ldots x_k)^2\prod_{i=1}^k \exp(-x_i^2/2)
$$
This is the density of the eigenvalues of a random $(k,k)$
matrix from the GUE, and classical results \cite{and}
naturally lead to the first part of the  conjecture. Proving it
would involve an exchange of limits, which does not
seem to have been achieved in the literature so far.\\

There is a link between rectangular Young tableaux and a particle system known 
as the TASEP, see \cite{romik}. In this view, a phenomenon of arctic circle 
arises when the rectangle is a square. In the general case, the arctic curve 
is no longer a circle but it is still algebraic. To compute an equation of
this curve, one has to determine the  parameters of the  Jacobi ensemble 
associated with the rectangle using Theorem 4, and then compute the 
values  $\lambda_\pm$ for this Jacobi ensemble as in \cite{dum}.

The horizontal strips of the TASEP diagram where this arctic curve 
appears correspond to diagonals of the Young tableau. Moreover, the places of the 
vertical steps inside a horizontal strip correspond to the integers
in the diagonal of the Young tableau, and Theorems 3 and 4 tell us that
these vertical steps are asymptotically distributed like the eigenvalues 
of a Jacobi ensemble.
A similar result has been established by Johansson and Nordenstam
for domino tilings of the Aztec diamond 
\cite{joh}, where the Jacobi ensemble is replaced by the GUE.

In the case of a GUE of size $n$, the fluctuations of the eigenvalues in the 
bulk are of order $\sqrt{\log n/n}$ and are asymptotically gaussian \cite{gust}.
It is not clear whether the same behaviour occurs 
in our context.

A model of corners of Jacobi ensembles was studied recently 
by Borodin and Gorin \cite{bor}, 
who showed a convergence to the gaussian free field.
Their model is slightly different from ours, but it would be interesting to 
know whether their results could be transposed in our case.

Finally, the method used here can be applied to generate at random a standard 
filling of a general polyomino: compute the conditional densities of 
the diagonals, which will be polynomials given by multiple integrals,
and then use a generating algorithm as in Section 2. Of course, the problem 
is that for a general polyomino, the corresponding polynomials will not have 
a simple form as for rectangles. 

{\bf Acknowledgements}
I thank Florent Benaych-Georges, C\'edric Boutillier, Benoit Collins
and Alain Rouault for useful references.


\begin{thebibliography}{999}

\bibitem{and}
Anderson, G., Guionnet, A., Zeitouni, O.,
{\em An introduction to random 
matrices.} Cambridge Studies in Advanced Mathematics, {\bf 118}. Cambridge 
University Press, Cambridge, 2010

\bibitem{bary}
Baryshnikov, Y., GUEs and queues. {\em Probab. Theory Related Fields}
{\bf 119} (2001), no. 2, 256-274.

\bibitem{bor}
Borodin, A., Gorin, G., 
General beta Jacobi corners process and the Gaussian Free Field,
{\em Communications on Pure and Applied Mathematics}
{\bf 68},  (2015),1774-1844.

\bibitem{coll}
Collins, B.,
Product of random projections, Jacobi ensembles and universality 
problems arising from free probability. 
{\em Probab. Theory Related Fields} {\bf 133} (2005), no. 3, 315-344. 

\bibitem{dum}
Dumitriu, I., Paquette, E.,
Global fluctuations for linear statistics of $\beta$-Jacobi ensembles.
{\em Random Matrices Theory Appl.}{\bf  1} (2012), no. 4

\bibitem{forr}
Forrester, P. J., {\em Log-gases and random matrices.} 
London Mathematical Society Monographs Series, {\bf 34}. 
Princeton University Press, Princeton, NJ, 2010.


\bibitem{gust}
Gustavsson, J.,
Gaussian fluctuations of eigenvalues in the GUE.
{\em Ann. I. H. Poincar\'e} (2005) 151-178.

\bibitem{joh}
 Johansson, K.; Nordenstam, E., 
Eigenvalues of GUE minors. {\em Electron. J. Probab.} {\bf 11} (2006), 
no. 50, 1342-1371. 

\bibitem{mar}
Marchal, P., 
Permutations with a prescribed descent set {\em Preprint}

\bibitem{muir}
Muirhead, R. J.,
{\em Aspects of multivariate statistical theory.}
Wiley Series in Probability and Mathematical Statistics. John Wiley and Sons, 
Inc., New York, 1982.

\bibitem{pitt}
Pittel, B., Romik, D.,
Limit shapes for random square Young tableaux.
{\em Adv. in Appl. Math.}{\bf 38} (2007), no. 2, 164-209. 

\bibitem{romik}
Romik, D., Arctic circles, domino tilings and square Young tableaux. 
{\em Ann. Probab.}{\bf 40} (2012), no. 2, 611-647.



\end{thebibliography}
\end{document}